\documentclass[12pt]{amsart}%{article}%
\usepackage{amsmath,amstext,amsfonts,amsthm,amssymb}
\usepackage{eucal}
\usepackage{diagrams}                                           %
\diagramstyle[scriptlabels,PostScript=dvips]               %
\newarrow{Dash}{}{dash}{}{dash}{}         
\renewcommand{\subsubsection}[1]{\addtocounter{subsubsection}{1}
{\ \\[3pt]\bf \thesubsubsection. \  #1} }
\swapnumbers

\newtheorem{THM}[subsection]{Theorem}

\newtheorem{PRP}[subsection]{Proposition}
\newtheorem{CRL}[subsection]{Corollary}

{  \theoremstyle{definition}

}
\newcommand{\Aut}{\operatorname{Aut}}
\newcommand{\Spec}{\operatorname{Spec}}
\newcommand{\bfR}{{\boldmath R}} %% R for right derived functors
\newcommand{\Hom}{\operatorname{Hom}}
\newcommand{\Char}{\operatorname{char}}
\newcommand{\cO}{\mathcal{O}}

\newcommand{\C}{\mathbb{C}}
\newcommand{\E}{\mathbb{E}}
\renewcommand{\P}{\mathbb{P}}

\newcommand{\Z}{\mathbb{Z}}

\begin{document}

\title[]{Subcanonical coordinate rings are Gorenstein}
\author{V.~Hinich}
\address{Department of Mathematics, University of Haifa,    
Mount Carmel, Haifa 31905,  Israel}    
\email{hinich@math.haifa.ac.il}    
\author{V.~Schechtman}
\address{Institut de Math\'ematiques, Universit\'e Paul Sabatier,
118 Route de Narbonne, 31062 Toulouse Cedex 9, France}
\email{schechtman@math.ups-tlse.fr}
%\begin{abstract}
%We prove that the coordinate rings of {\it subcanonical} smooth projective 
%varieties are Gorenstein and have rational singularities. 
%Combined with a classical theorem of Avramov and Golod, this enables to 
%strenghten and extend some recent results of Gorodentsev, 
%Khoroshkin and Rudakov. 
%\end{abstract}
\maketitle    

%
%  WRITE HERE
%
\begin{quote}
{\em To our teacher Evgeny Solomonovich Golod, with gratitude} 
\end{quote}

\vspace{1cm}

%\vskip{1cm}

\begin{quotation}
{\em{\small In all examples we consider, {\em [the coordinate ring 
of $X$] }%\\
 is a Gorenstein ring; this property is one of the most powerful %\\
general tools we have in studying $X$ and its deformations.%\\
 It seems to us that this point is not adequatly appreciated.}}
 
\vspace{0.3cm} 

A. Corti, M. Reid, Weighted Grassmanians.  
\end{quotation}

\vspace{0.3cm}

\section{Introduction}

Let $i: X \rInto \P = \P(V)$ be a smooth connected projective variety 
embedded 
into a projective space (we are working over a fixed ground field $k$).  
Set $\cO_X(1) = i^*\cO_\P(1)$ and consider the coordinate algebra 
$$
A %= \bigoplus_{n=0}^\infty\ A_n 
= \bigoplus_{n=0}^\infty H^0(X,\cO_X(n)).
$$
By construction $A$ is identified with a quotient algebra $A = S/I$ where 
$S = Sym(V^*) = k[x_0,\ldots, x_{n-1}]$. 
The {\it Koszul homology} algebra is defined as 
$$
H(A) = \bigoplus_{p=0}^n\ Tor_p^S(A,k).
$$
This is a (bi)graded commutative $k$-algebra, finite dimensional as a 
$k$-vector space. 

In an inspiring paper \cite{GKR} Gorodentsev, Khoroshkin and Rudakov 
prove (among others) the following elegant result. Denote by $K_X$ the
canonical class of $X$.

\begin{THM}[(see~\cite{GKR}, Sect. 2)]
\label{thm:GKR}
Suppose that 
\begin{itemize}
\item[(a)] there exists a natural $N$ such that $K_X = \cO_X(-N)$; 
\item[(b)] $H^i(X,\cO_X(n)) = 0$ for all $n \in \Z$ and $0 < i < d := \dim X$. 
\end{itemize}
Then $H(A)$ is Frobenius. 
\end{THM}

Here {\it Frobenius} means that there exists a nondegenerate bilinear pairing 
\newline $\langle\ ,\ \rangle:\ H(A) \times H(A) \rTo k$, 
suitably compatible with the gradings, such that 
$\langle ab, c\rangle = \langle a, bc\rangle$. 

The proof in {\it op. cit.} is very nice; it uses the "sphericity" of certain 
spectral sequence. 

In this note we would like to look at this result from a slightly different 
perspective.  
Our point of departure is a fundamental result by Avramov and Golod,

\cite{AG}:

\begin{THM}
 $H(A)$ is Frobenius if and only if $A$ is Gorenstein. 
\end{THM}
In fact, Avramov and Golod work in the local situation; the passage to our 
graded context presents no difficulties. Indeed, according to {\it op. cit.}, 
$H(A)$ is Frobenius iff the localisation of $A$ at $0$ is Gorenstein; however, 
$A$ is smooth outside this ideal, so this is equivalent to $A$ being 
Gorenstein.  

So our question reduces to the
Gorenstein property of $A$. 

Let us say, following \cite{GKR}, 
that $X \subset \P$ is {\it subcanonical} if the condition (a) of 
Theorem~\ref{thm:GKR} is satisfied. In the present note we prove the following 

\begin{THM}\label{thm:1} 
Assume $\Char(k)=0$.
If $X\subset \P$ is subcanonical then $A$ is Gorenstein and has rational 
singularities. 
\end{THM}

We establish this using certain {\it Key Lemma} from \cite{H} (see Proposition 
\ref{keyl}) giving a
sufficient condition for a singularity being Gorenstein and rational. 
The proof of this lemma uses Grauert-Riemenschneider theorem, 
and hence the characteristic zero assumption. (On the contrary, although 
Gorodentsev et al. 
assume $k=\C$, their proof of \ref{thm:GKR} works over an arbirtary field). 

\begin{CRL}\label{crl:1}
If $X\subset \P$ is subcanonical then $H(A)$ is Frobenius. 
\end{CRL}

So, the condition (b) of Theorem \ref{thm:GKR} is superfluous if 
$\Char k = 0$.  

The main objects of study in {\it op cit.} are 
{\it highest weight orbits} of a semisimple algebraic 
group $G$. For such $X$ the authors of \cite{GKR} prove that (b) follows from 
(a). 

In this case we prove that subcanonicity is {\it equivalent} to 
the  Gorenstein property of $A$:   

\begin{THM}\label{thm:2}
Let $ X \subset \P(V)$ be the 
projectivisation of the highest weight orbit in an irreducible finite 
dimensional representation $V$ of a semisimple group $G$. 
This embedding is subcanonical if and only if the corresponding coordinate 
ring $A$ is Gorenstein (so, iff $H(A)$ is Frobenius).  
\end{THM}

\subsection{Acknowledgement}
This note was written during a visit of the first author to the 
{\it Institut de Math\'ematiques de Toulouse}. He thanks this Institute for 
the hospitality.     

\newpage

\section{Proof of Theorem \ref{thm:1}}

We keep the notation of the Introduction. 
The affine variety $Z:=\Spec(A)$ is the cone over $X$; 
therefore it is nonsingular outside $0$. 
It has a very nice desingularization $Y$ which is
the total space of the vector bundle $\E=\cO_X(-1)$. Let
\begin{equation}\label{eq:p}
p:Y=\Spec(Sym_{\cO_X}(\E^*))\rTo X
\end{equation}
be the projection. 

The embedding $\cO_{\P(V)}(1)\rTo V$ defines an embedding
$Y\rTo X\times V$; the projection to the second factor has image
$Z\subset \Spec (Sym\ V^*) = V$ and the map
\begin{equation}\label{eq:pi}
 \pi:Y\rTo Z
\end{equation}
is a desingularization.

Recall the following 

\begin{PRP}[see~\cite{H}]\label{keyl}
Let $\pi:Y\rTo Z$ be a proper birational map with $Y$ smooth and $Z$ normal.
Let $\omega_Y$ be the sheaf of higher differentials on $Y$. Assume there exists
a morphism $\phi:\cO_Y\to\omega_Y$ such that 
$\pi_*\phi:\pi_*\cO_Y\to\pi_*\omega_Y$ is an isomorphism. Then $Z$ is Gorenstein
and has rational singularities.
\end{PRP}

We wish to apply this to our desingularization $\pi:Y\to Z$. Note that
$Z=\Spec(A)$ is normal.

The short exact sequence of vector bundles on $Y$
\begin{equation}\label{eq:ses}
0\rTo p^*\E\rTo T_Y\rTo p^*T_X\rTo 0
\end{equation}
yields an isomorphism 
\begin{equation}\label{eq:omegay}
\omega_Y=p^*(\omega_X\otimes\E^*).
\end{equation}
We wish to calculate the global sections of $\omega_Y$. First of all,
we have
\begin{equation}\label{eq:p-omegay}
p_*\omega_Y=p_*p^*(\omega_X\otimes\E^*)=\omega_X\otimes\E^*\otimes 
Sym_{\cO_X}\E^*=\bigoplus_{n\geq 1}\omega_X\otimes\cO_X(n)
\end{equation}
since $p$ is an affine morphism. 

\subsection{Proof of Theorem \ref{thm:1}} 

Let $\omega_X=\cO_X(-N)$.
One has an obvious map
$$ \cO_X=\omega_X\otimes\cO_X(N)\rInto
\bigoplus_{n\geq 1}\omega_X\otimes\cO_X(n)=p_*\omega_Y$$
which gives by adjunction a map $\phi:\cO_Y\rTo\omega_Y$. 

We will check now that $\phi$ induces an isomorphism of the global sections.
Applying to $\phi$ the direct image functor $p_*$
we get a morphism
\begin{equation}\label{eq:p-phi}
p_*(\phi): \bigoplus_{n\geq 0}\cO_X(n)\rTo 
\bigoplus_{n\geq 1}\omega_X\otimes\cO_X(n)
\end{equation}
which is obviously a map of $p_*(\cO_Y)$-modules. By definition it carries
$1\in p_*(\cO_Y)$ to a generator of $\omega_X(N)=\cO_X$, so the map 
$p_*(\phi)$
carries isomorphically the summand $\cO_X(n)$ of the left-hand side to 
the summand $\omega_X\otimes\cO_X(N+n)$ of the right-hand side.
For $n<N$ one has on the right-hand side
$$ \Gamma(X,\omega_X\otimes\cO_X(n))=\Gamma(X,\cO_X(n-N))=0,$$
so $p_*(\phi)$ induces an isomorphism of the global sections.

\section{Homogeneous case}

Let now $G$ be a semisimple Lie group, $V$ a simple finite dimensional
highest weight $G$-module, $v\in V$ be a highest weight vector.
Let $P$ be the stabilizer of $\C v$ in $\P(V)$. This is a parabolic
subgroup of $G$. A $G$-equivariant embedding $i:X:=G/P\rTo\P(V)$ is 
induced.

The closure $Z$ of $Gv$ is a cone in $V$.  We have $Z=\Spec(A)$
where $A$ is the homogeneous coordinate ring of $X=G/P$ with respect to $i$.

In this case the converse of the theorem~\ref{thm:1} is valid. One has

\begin{THM}
The space $Z$ is Gorenstein iff $\omega_X=\cO_X(-N)$ for some $N$.
\end{THM}

Note that the conclusion of the Theorem is not true for an arbitrary 
(nonhomogeneous) $X$ (for example it follows easily from the results of Mukai 
\cite{M}  
that a generic curve of genus $7$ embedded canonically in $\P^6$ has a 
Gorenstein coordinate ring).

\begin{proof}

The dualizing complex of $Z$ can be calculated as
\begin{equation}\label{omega-Z}
\omega_Z=\bfR\Hom_{SV^*}(A,SV^*)[\dim V-\dim Z]
\end{equation}
(the shift is chosen so that $\omega_Z$ is concentrated in degree $0$ when $A$
is Cohen-Macaulay).

Its cohomology keeps the grading of $SV^*$ and $A$; therefore, if $A$
is Gorenstein so that $\omega_Z$ is an invertible $A$-module, it has
to be isomorphic to $A$.

Choose an isomorphism $\theta:A\rTo\omega_Z$.

We now apply the Duality isomorphism, see \cite{Ha}, VII.3.4, to the proper
morphism $\pi:Y\to Z$. It gives, in particular, an isomorphism
\begin{equation}\label{eq:duality}
\Hom_{D(Y)}(F,\pi^!G))\rTo^\sim\Hom_{D(Z)}(\bfR\pi_*F,G)
\end{equation} 
for any $F\in D^-_{qc}(Y),\ G\in D^+_c(Z)$.

We apply this to $F=\cO_Y$ and $G=\omega_Z$.
By a general result of Kempf~\cite{K} $Z$ has rational singularities,
so $\bfR\Gamma(Y,\cO_Y)=\Gamma(Y,\cO_Y)=A$. Moreover, 
$\pi^!(\omega_Z)=\omega_Y$. Thus, Duality isomorphism gives us
\begin{equation}\label{eq:ourduality}
\Hom_{D(Y)}(\cO_Y,\omega_Y))\rTo^\sim\Hom_{D(Z)}(\cO_Z,\omega_Z).
\end{equation}

We see that the map $\theta:A\to\omega_Z$ is adjoint to a map
$\theta_Y:\cO_Y\to\omega_Y$ which in turn can be rewritten as a morphism
\begin{equation}\label{eq:theta-x}
\theta_X:\cO_X\to p_*(\omega_Y)=\bigoplus_{n\geq 1}\omega_X(n).
\end{equation}

We intend to prove now that each direct component 
$\theta_{X,n}:\cO_X\to\omega_X(n)$
is either isomorphism or vanishes. This will immediately imply the theorem.

Note that the formula~(\ref{omega-Z}) shows that the group $G$ naturally 
acts on $\omega_Z$.  We claim that $\theta:A\to\omega_Z$ is necessarily
$G$-equivariant.

In fact, the $G$-action on $A$-module $\omega_Z$ is compatible with $G$-action
on $A$:
$$ g(ax)=g(a)g(x),\ g\in G,\ a\in A,\ x\in \omega_Z.$$
Another $G$-module structure on $\omega_Z$ compatible with the $G$-action on 
$A$ is given by $\theta$. These two actions define two group homomorphisms
$$ \rho_1,\rho_2:G\rTo\Aut_\C(\omega_Z).$$
The ``difference'' between the two defined by the formula
$$\rho_{12}:g\mapsto \rho_1(g^{-1})\circ\rho_2(g)$$
gives rise to a crossed homomorphism $\rho_{12}:G\to\Aut_A(\omega_Z)=\C^*.$
Since the action of $G$ on $\C^*$ is trivial and $G$ is semisimple,
$\rho_{12}$ is trivial, which means that $\theta$ is $G$-equivariant.

Let us show that the maps $\theta_Y$ and $\theta_X$ obtained from $\theta$
via Duality isomorphism, are also $G$-equivariant.

Choose $g\in G$ and let $g_X:X\to X,\ g_Y:Y\to Y,\ g_Z:Z\to Z$ denote the 
corresponding automorphisms of the varieties.

An action of $g\in G$ on $\cO_Z$ and $\omega_Z$ are expressed as isomorphisms
$g_Z^*(\cO_Z)\to\cO_Z$ and $g_Z^*(\omega_Z)\to\omega_Z$. Since $\theta$ is 
equivariant, it gives rise to a commutative diagram
\begin{equation}
\begin{diagram}
g_Z^*(\cO_Z) & \rTo^{g_Z^*\theta} & g_Z^*(\omega_Z) \\
\dTo & & \dTo \\
\cO_Z & \rTo^\theta & \omega_Z
\end{diagram}
\end{equation}

The map $\theta_Y$ can be described as the composition
$$ \cO_Y\rTo\pi^!\bfR\pi_*(\cO_Y)=\pi^!\cO_Z\rTo\pi^!\omega_Z,$$
so that it suffuces to check that the first morphism is $G$-equivariant.
The latter can be expressed as the commutativity of the diagram
\begin{equation}
\begin{diagram}
g_Y^*(\cO_Y) & \rTo & g_Y^*(\pi^!\bfR\pi_*(\cO_Y)) \\
\dTo & & \dTo \\
\cO_Y & \rTo & \pi^!\bfR\pi_*(\cO_Y)
\end{diagram}
\end{equation}
for each $g\in G$, and this follows from the relations
$$ g^*_Y\pi^!=\pi^!g^*_Z,\quad g^*_Z\bfR\pi_*=\bfR\pi_*g^*_Y.$$

All this proves that $\theta_Y$ is $G$-equivariant; the similar fact
for $\theta_X$ is even more transparent.

We have already understood that the components $\theta_{X,n}$
of the map \\
$\theta_X:\cO_X\rTo\bigoplus \omega_X(n)$
are $G$-equivariant. This implies that 
the map of fibers at $1P\in G/P$ is $P$-equivariant.
The fibers are one-dimensional representations of $P$; any $P$-morphism
is either zero or an isomorphism. This proves the theorem.

\end{proof}

\end{document}